\documentclass[a4paper,12pt]{amsart}
\usepackage{amssymb,amsthm}
\usepackage[all]{xy}

\begin{document}

\newcommand{\opp}{\bowtie }
\newcommand{\po}{\text {\rm pos}}
\newcommand{\supp}{\text {\rm supp}}
\newcommand{\End}{\text {\rm End}}
\newcommand{\diag}{\text {\rm diag}}
\newcommand{\Lie}{\text {\rm Lie}}
\newcommand{\Ad}{\text {\rm Ad}}
\newcommand{\car}{\mathcal R}
\newcommand{\Tr}{\rm Tr}
\newcommand{\Spec}{\text{\rm Spec}}

\def\ge{\geqslant}
\def\le{\leqslant}
\def\a{\alpha}
\def\b{\beta}
\def\c{\chi}
\def\g{\gamma}
\def\G{\Gamma}
\def\d{\delta}
\def\D{\Delta}
\def\L{\Lambda}
\def\e{\epsilon}
\def\et{\eta}
\def\io{\iota}
\def\o{\omega}
\def\p{\pi}
\def\ph{\phi}
\def\ps{\psi}
\def\r{\rho}
\def\s{\sigma}
\def\t{\tau}
\def\th{\theta}
\def\k{\kappa}
\def\l{\lambda}
\def\z{\zeta}
\def\v{\vartheta}
\def\va{\varphi}
\def\x{\xi}
\def\i{^{-1}}

\def\mapright#1{\smash{\mathop{\longrightarrow}\limits^{#1}}}
\def\mapleft#1{\smash{\mathop{\longleftarrow}\limits^{#1}}}
\def\mapdown#1{\Big\downarrow\rlap{$\vcenter{\text{$\scriptstyle#1$}}$}}
\def\mapup#1{\Big\uparrow\rlap{$\vcenter{\text{$\scriptstyle#1$}}$}}

\def\ca{\mathcal A}
\def\cb{\mathcal B}
\def\cc{\mathcal C}
\def\cd{\mathcal D}
\def\ce{\mathcal E}
\def\cf{\mathcal F}
\def\cg{\mathcal G}
\def\ch{\mathcal H}
\def\ci{\mathcal I}
\def\cj{\mathcal J}
\def\ck{\mathcal K}
\def\cl{\mathcal L}
\def\cm{\mathcal M}
\def\cn{\mathcal N}
\def\co{\mathcal O}
\def\cp{\mathcal P}
\def\cq{\mathcal Q}
\def\car{\mathcal R}
\def\cs{\mathcal S}
\def\ct{\mathcal T}
\def\cu{\mathcal U}
\def\cv{\mathcal V}
\def\cw{\mathcal W}
\def\cz{\mathcal Z}
\def\cx{\mathcal X}
\def\cy{\mathcal Y}

\def\ccg{\mathfrak g}
\def\ccb{\mathfrak b}
\def\cch{\mathfrak h}
\def\cck{\mathfrak k}
\def\ccp{\mathfrak p}
\def\ccn{\mathfrak n}
\def\cct{\mathfrak t}
\def\ccl{\mathfrak l}
\def\cci{\mathfrak i}
\def\ccq{\mathfrak q}
\def\ccu{\mathfrak u}
\def\ccs{\mathfrak s}
\def\cco{\mathfrak o}
\def\ccd{\mathfrak d}
\def\cca{\mathfrak A}

\def\CR{\mathbb R}

\def\tz{\tilde Z}
\def\tl{\tilde L}
\def\tc{\tilde C}
\def\ta{\tilde A}
\def\tx{\tilde X}

\newtheorem{theorem}{Theorem}[section]
\newtheorem{lem}[theorem]{Lemma}
\newtheorem{cor}[theorem]{Corollary}
\newtheorem{prop}[theorem]{Proposition}
\newtheorem{thm}[theorem]{Theorem}
\newtheorem*{rmk}{Remark}
\newtheorem{eg}[theorem]{Example}
\newtheorem{defi}[theorem]{Definition}
\newtheorem{conj}[theorem]{Conjecture}
\newtheorem{example}[theorem]{Example}
\newtheorem{ack}{Acknowledgement}

\newenvironment{thmref}{\thmrefer}{}
\newcommand{\thmrefer}[1]{\renewcommand\thetheorem
{\protect\ref{#1}}\addtocounter{theorem}{-1}}

\author{Chuying Fang}
\address{Department of Mathematics, Massachusetts Institute of Technology, Cambridge, MA 02139, USA}
\email{cyfang@alum.mit.edu}
\title[Ad-nilpotent Ideals of Minimal Dimension]{Ad-nilpotent Ideals of Minimal Dimension}

\begin{abstract}
We use Jacobson-Morozov theorem to prove Sommers' conjecture about  the lower bounds of ad-nilpotent ideals with the same associated orbit. More precisely,  for each nilpotent orbit, we  construct some minimal ad-nilpotent ideals corresponding to their associated orbit.

For classical groups of  type $A_{n-1}$, we  get an explicit formula for the minimal dimension based on partitions of $n$.
\end{abstract}

\maketitle

\section{Introduction}
Let $G$ be a complex simple Lie group with Lie algebra  $\ccg$. Fix a Borel subgroup of $G$.  Let $\ccb$ be the Lie algebra of $B$ and $\ccn$ the nilradical of $\ccb$.

An ideal of $\ccb$ is called {\it ad-nilpotent}  if it is contained in the nilradical $\ccn$ (sometimes it is called a $B$-stable ideal) . Ad-nilpotent ideals have many applications in the study of affine Weyl groups, hyperplane arrangements, sign types, and nilpotent orbits. We refer to  \cite{CP1} \cite{GS} \cite{Mi}\cite{Pa1} \cite{So1}\cite{So2} for recent result about ad-nilpotent ideals and their applications. In particular, ad-nilpotent ideal was a main tool of  Mizuno \cite{Mi} to study  the conjugate classes of nilpotent elements for exceptional groups. Some further results in this direction were obtained by Kawanaka \cite{Kaw1},   Gunnells-Sommers \cite{GS} and Sommers \cite{So2}.  

Suppose $I$ is  an ad-nilpotent ideal. Consider the map $G\times_B I \rightarrow \ccg$, which is a restriction of the moment map $G\times_B \ccn\rightarrow \ccg$. Its image is the closure of one unique nilpotent orbit and we denote it by $\co_I$. This orbit is called the associated orbit of the ideal $I$.   This induces a map  from the set of ad-nilpotent ideals to the set of nilpotent orbits.  By Jacobson-Morozov theorem, for each nilpotent orbit $\co$, there always exists an ad-nilpotent ideal  associated to $\co$. Namely, the map is a surjection.


Sommers showed in  \cite{So2} that   the dimensions of the ideals with  the same associated orbit $\co$ have a lower bound $m_\co$.
\begin{prop} \cite[Prop5.1]{So2}
Let $\co$ be a nilpotent orbit. Let $I$ be an ad-nilpotent ideal whose associated orbit is $\co$ and $X$ be an element that lies both in the ideal $I$ and $\co$. Then $$\dim I \ge \dim B- \dim  B_{G_X}$$ where $G_X$ is the centralizer  of $X$ in $G$ and $B_{G_X}$ is a Borel subgroup of $G_X$.
\end{prop}

Notice that $\dim B_{G_X}$ is independent of the choice of $X$, but only depends on the orbit $\co$.  Set  $m_\co=\dim B - \dim B_{G_X} $.  Then Sommers conjectured that

\begin{conj}\label{conj}
For each nilpotent orbit $\co$, there exists an ideal $I$ with $\co_I= \co $ and  $\dim I= m_\co$.
\end{conj}

The main purpose of this paper is to study the minimal dimension of ad-nilpotent ideals with the same associated orbit and prove the conjecture for classical groups.  For exceptional groups, this conjecture was implicitly proved by the  work of Kawanaka  \cite{Kaw1} and Mizuno \cite{Mi}. 

We give a brief  outline of  this paper. In section 2, we introduce the standard triples for the nilpotent orbits and construct Dynkin ideals from the Dynkin element of the orbit. In most cases, the Dynkin ideals  are not ideals of minimal dimension, but can provide some useful information about the minimal ideals. In sections 3, we  give a  general strategy to construct the ideals of minimal dimension and construct some ideals of the  minimal dimension in the case of type $A$.    Section 4 is an application of section 3, where we give a formula for the dimension of minimal ideals for the type $A_{n-1}$  based on partitions of $n$ and then derive that if $\co_1$ and $\co_2$ are nilpotent orbits with $\co_1$ is contained in the closure of $\co_2$, then $m_{\co_1}\le m_{\co_2}$. In sections 5, 6, and 7, the explicit construction of minimal ideals for type $B$, $C$ and $D$ is given.

I would like to  express my deep gratitude to  my advisor David Vogan for his guidance, warm encouragement and many useful suggestions while  preparing the paper.  I would like to thank George Lusztig for his enjoyable lectures.  I  would also like to  thank Eric Sommers for both the email correspondence and conversations. 

\section{Notation and Preliminaries}

If $V$ is a subspace of $\ccg$ that's invariant under the action of $\cch$, we denote $\D(V)=\{ \a \in \D \mid \ccg_{\a} \in V\}$. If $S$ is a finite set, we denote by $|S|$ the cardinality of $S$.  For any $k \in \mathbb R$, let  $ \lfloor k\rfloor $ be the  largest integer less than or equal to $k$ and let $ \lceil k\rceil$ be the smallest integer not less than $k$.

Notice that any ad-nilpotent ideal is completely determined by its underline set of roots. For any ad-nilpotent ideal $I$, let

Before we come to the proof of the conjecture of Sommers, let's first recall some results about standard triples.

Let $\{H, X, Y\}$ be a $\bold {standard}$ triple (see  \cite{CM}) of $\ccg$, satisfying:
$$[H,X]=2X, \qquad [H,Y]=-2Y, \qquad[X,Y]=H.$$
We call $X$ (resp. $Y$) the {\it nilpositive } (resp. {\it nilnegative }) element and $H$ the {\it characteristic} of the triple $\{H,X,Y\}$. In particular,   after conjugation by some  element of $G$,  we can assume that  $H \in \cch$ and  $H$ is dominant, i.e. $\a(H) \ge 0 $, for all $\a \in \D^+$. Such $H$ is uniquely determined by the nilpotent orbit $\co_X$ and is called the {\it Dynkin element} for  $\co_X$.   There is an $H$-eigenspace decomposition of $\ccg$: $$\ccg=\oplus_{i\in \mathbb Z}\ccg_{H,i}, \text{ where } \ccg_{H,i}= \{Z \in \ccg \mid [H, Z]= iZ\}, i \in \mathbb Z.$$

It's also shown in  \cite{So2} that \[\tag{2.1.1} \dim B_{G_X} =\dim  \ccg_{H,1}+ \frac{1}{2}[\dim (\ccg_{H,0})+\dim (\ccg_{H,2})+ \text{rank } G_X]. \]

Let $\ccq_{H,i}=\oplus_{j \ge i}\ccg_{H,j}$.  Then, $\ccq_{H,i}$ is an ad-nilpotent ideal and  $X  \in \ccq_{H,2}$.  We call the ideal $\ccq_{H,2}$ the {\it Dynkin ideal} for the orbit $\co_X$. We may write  $X$ as  $X_{\a_1}+ X_{a_2}+ \dots +X_{\a_k}$, where $X_{\a_i}$ is a root vector and $\a_i(H)=2$. Since $H$ is dominant, each $\a_i$ is a positive root.

Consider the adjoint action $ad_X: \ccg_{H,0} \rightarrow \ccg_{H, 2}$, which sends any $Z \in \ccg_{H,0} $ to $[X, Z] \in \ccg_{H,2}$. Since $\ccg_{H,0}$ is a Levi subalgebra of $\ccg$ containing $\cch$, there is a direct sum decomposition:  $\ccg_{H,0}=\cch \oplus \ccg_0^+ \oplus \ccg_0^-  $, where $\ccg_0^+= \ccn \cap \ccg_{H,0} $ and $\ccg^-_0=  \ccn^- \cap \ccg_{H,0}$.   It's obvious that $ad_X(\ccg_{H,0})=ad_X(\cch)+ ad_X(\ccg^+_0)+ad_X(\ccg^-_0)$.   If we impose some additional restrictions on the set $\{\a_1, \a_2, \dots \a_k\}$,  we will have a direct sum decomposition of $\ccg_{H,2}$. First we introduce the notion of antichain.

\begin{defi} \cite{St1}
An antichain  $\G$ of the root poset $(\D^+,<)$ is a set of pairwise incomparable elements, i.e: for any $\a, \b \in \G$, $\a- \b \notin Q^+$, where $Q^+= \{\sum_{i=1}^n n_i\a_i \mid n_i \in \mathbb N\}$ is the positive part of the root lattice.
\end{defi}

By the definition of generators of an ad-nilpotent ideal, a set $\G=\{\g_1, \dots, \g_l\}$ is a set of generators of some ad-nilpotent ideal if and only if $\g_i -\g_j \notin Q^+$. Therefore, the set of generators for ad-nilpotent ideals is in bijection with the set of antichains of the root poset.

We need the following result of Kostant. (see  \cite{Kos}):

\begin{thm}\label{Kos}
Let $Q_{H,2}$, $G_{H,0}$ be the  closed  connected Lie subgroup of $G$ with the Lie algebras $\ccq_{H,2}, \ccg_{H,0}$ respectively. Let $\co_X$ be the G-orbit of $X$ and $ o_X$ the $G_{H,0}$-orbit of $X$. Then

(1)  $o_X$ is open, dense in $\ccg_{H,2}$.

(2) $o_X=\co_X \cap \ccg_{H,2}$.

(3) $(Q_{H,2}G_{H, 0}) \cdot X=\co_X \cap \ccq_{H,2}=o_X +\ccq_{H,3}$. In particular, $(Q_{H,2}G_{H, 0}) \cdot X$ is open and dense in $\ccq_{H,2}$.
\end{thm}

As a consequence,  the Dynkin ideal $\ccq_{H,2}$ has associated orbit $\co_X$.

\begin{lem} \cite[Prop2.10]{Pa1}\label{pa}
Let $\G$ be a subset of $\D^+$. If for any roots $\a ,\b \in \G$, $\a-\b \notin \D$, then the elements of $\G$ are
linearly independent and hence $ |\G| \le dim(\cch)$.
\end{lem}
Proof. Since $\a-\b \notin \D$, $(\a, \b)\le 0$. This means that the angle between  any pair of roots in   $\G$ is non-acute. Since all the roots  in $\G$ lie in the same open half-space of $V$, they are linearly independent. \qed

\begin{rmk} If $\G$ is an antichain, then $\G$ satisfies the assumption of lemma \ref{pa}, hence elements in an antichain are linearly independent. In Panyushev's original statement, he assumed that $\G$ is an antichain. But from his proof, the weaker condition that $\a- \b \notin \D$ is sufficient for this lemma. In type $A_{n-1}$ and $D_n$, we can get an antichain but in type $C_n$ case,   we can only get a set $\G$ satisfying  the weaker condition of Lemma \ref{pa}.\end{rmk}

\begin{lem}\label{trip}
Suppose $\G$ is a subset of $\D^+$ as in Lemma \ref{pa}. Let $\{H_\a, X_\a, Y_\a\}$ be a standard triple that corresponds to  $\a \in \G$. Suppose $H$ lies in the span of all $\{H_\a\}_{\a \in \G}$ so that $\a(H)=2$ for all $\a \in \G$.  Let $X=\sum_{\a \in \G}X_\a$. There exists an element $Y$, such that  $\{H, X, Y \}$  is a standard triple. 
\end{lem}
\begin{rmk}The main idea of the proof comes from   \cite[4.1.6]{CM}. \end{rmk}

Proof. Since $\G$ is a subset as in Lemma \ref{pa}, $\{H_\a \mid \a \in \G\}$ is linearly independent. Also $H$ lies in the subspace of $\cch$ that's spanned by all $\{H_\a\}_{\a \in \G}$, therefore we may write $H$ as   $H=\sum_{\a \in \cc}a_\a H_\a$ for some $a_\a \in \mathbb R$. Let $Y= \sum_{\a \in \G} a_\a Y_\a$. Then $[H, X]=[H, \sum_{\a\in \G} X_\a]=\sum_{\a \in \cc} \a(H)X_\a=2X$. The last equality follows from $(2.2.1)$.  Similarly, $[H, Y]=-2Y$. Also $$[X,Y]= \sum_{\a \in \G}\sum_{\b \in \G}a_{\b}[X_\a, Y_\b]= \sum_{\b \in \G}a_\b[X_\b, Y_\b]=\sum_{\b \in \cc}a_{\b}H_{\b}=H.$$ The  second equality comes from the fact that  $\a-\b \notin \D$ and $[X_\a, Y_\b]=0$ for $\a \neq \b$.\qed

\

\begin{prop}
Suppose that the set $\{\a_1, \a_2, \dots \a_k\}$ is an antichain of  $(\D^+, <)$.  Then $\ccg_{H,2}= ad_X(\ccg_{H,0}) $ and $\ccg_{H,2}= ad_X(\cch) \oplus ad_X(\ccg^+_0) \oplus ad_X(\ccg^-_0)$. Moreover, $ad_{X}(\cch)$, $ad_X(\ccg_0^+)$ and $ad_X(\cch)$ are invariant under the adjoint  action of $\cch$ and each can be written as a direct sum of root spaces.
\end{prop}

Proof.  By Theorem \ref{Kos}, the image of $ad_X$ is the whole space $\ccg_{H,2}$. What remains to prove is that it is a direct sum decomposition. Suppose there exist three elements $H_1 \in \cch$ , $Z \in \ccg_0^+$, $U \in \ccg^-_0$ and $[X, H_1]+[X, U]+[X, Z]=0$.

Suppose $\a_i(H_1) \neq 0$ for  some $\a_i$. Since $[X, H_1]=\sum_{i=1}^k -\a_i(H_1)X_{\a_i} $, there is a nonzero summand in $[U, X]$ or $[Z, X]$ that lies in $\ccg_{\a_i}$. Without loss of generality, we may assume $[U_{\b}, X_{\a_j}] \in \ccg_{\a_i}$, where $U_{\b}$ is a summand of $U$ and $X_{\a_j}$ is a summand of $X$. Then $\a_i= \a_j + \b$, which implies that $\a_i< \a_j$ and contradicts the assumption that $\a_i$ and $\a_j$ are incomparable.

Otherwise $\a_i(H)=0$ for $1 \le i \le k$.  By similar argument, we can prove that $ad_{X}(\ccg_0^+) \cap ad_X(\ccg_0^- )=0 $, which shows that  $\ccg_{H,2}= ad_X(\cch) \oplus ad_X(\ccg^+_0) \oplus ad_X(\ccg^-_0)$.

Suppose that $X_\a$ is a root vector in $\ccg_{H,2}$. By lemma \ref{pa}, $\{\a_1, \a_2, \dots, \a_k\}$ are linearly independent. If $\a \in \{\a_1, \a_2, \dots, \a_k\}$, then $\a$ lies in $ad_X(\cch)$. If $\a=\a_i +\b$, where $\b \in \D(\ccg_0^+) $, from the proof above, $X_\a$ can not appear in the summands of $ad_X(\cch)$ and $ad_X(\ccg_0^-)$. Similarly, elements in $ad_X(\ccg_0^-)$ has the form $\sum X_{\a_i+ \b_i}$, where $\b_i \in \D(\ccg^-_0)$. Therefore, $X_\a \in ad_X(\ccg_0^+)$ if $\a=\a_i+\b$ and $\b \in \D(\ccg_0^+)$ and $X_\a \in ad_X(\ccg_0)^-$ if $\a=\a_i+\b$ and $\b \in \D(\ccg_0^-)$. This shows that $ad_X(\ccg_0^-)$, $ad_X(\ccg_0^+)$ and $ad_X(\cch)$ are $\cch$-invariant, which completes the proof.   \qed

\begin{defi}  A $\bold {partition}$ $\l$ is a  sequence of positive integers $$\l=[\l_1, \l_2, \dots, \l_p], \text{ where } \l_1 \ge \l_2 \dots \ge \l_p >0.$$ Each $\l_i$  is a part of  $\l$.  If $\l_1+\dots +\l_p=n$,  we say $\l$ is a partition of $n$ and write $\l \vdash n$.
\end{defi}


We recall the  parametrization of nilpotent orbits  in classical groups.

\begin{thm} \cite[5.1]{CM}\label{cm1}
(1)(Type $A_{n-1}$) Nilpotent orbits in $\ccs \ccl_n$ are in
one-to-one correspondence with the set $P(n)$ of partitions of $n$.

(2)(Type $B_n$) Nilpotent orbits in  $\ccs \cco_{2n+1}$ are in one-to-one correspondence with the set $P_1(2n+1)$ of partitions of $2n+1$ in which even parts occur with even multiplicity.

(3)(Type $C_n$) Nilpotent orbits in $\ccs \ccp_{2n}$ are in one-to-one correspondence with the set $P_{-1}(2n)$of partitions of $2n$ in which odd parts occur with even multiplicity.

(4)(Type $D_n$) Nilpotent orbits in $\ccs \cco_{2n}$ are in
one-to-one correspondence with the set $P_1(2n)$ of partitions of
$2n$ in which even parts occur with even multiplicity, except that
(`very even') partitions (those with only  even parts; each
having even multiplicity) correspond to two orbits.
\end{thm}

For each partition $\l$, we denote by $\co_\l$ the nilpotent orbit that corresponds to  $\l$ except the very even case in type $D_n$, in which we denote the two orbits by $\co^I_{\l}$ and $\co^{II}_\l$.

We will prove conjecture \ref{conj} by  constructing explicit minimal ideals
in the classical groups. In addition, we have an explicit formula
for the dimension of the minimal ideals in terms of partition.

Let's briefly discuss the main idea to construct minimal ideals. First let's recall the method to compute the weighted Dynkin diagram of a nilpotent orbit   in  \cite{CM}.  Given a partition $\l=[\l_1, \dots, \l_p]$ as a nilpotent orbit in type $A_{n-1}$, for each part $\l_i$, we take the set of integers $\{\l_i-1, \dots,1- \l_i\}$. Then we take the union of these sets and write it into a sequence $(h_1, h_2, \dots, h_n)$, where $h_1 \ge \dots \ge h_n$.  We assign the value $h_i-h_{i+1}$ to the $i$-th node of the Dynkin diagram of $A_{n-1}$. This gives us the weighted Dynkin diagram corresponding to $\l$.

Coming back to the construction of ad-nilpotent ideals, for each $h_i$ in the sequence $(h_1, \dots, h_n)$, we have to specify which part of the partition it comes from.  For example, if $h_{i_1}, \dots, h_{i_k}$ come from the same part of $\l$ in a descending order, then we may pick the roots $\{e_{i_j}-e_{i_{j+1}}\}_{j=1}^{k-1}$ to be  generators of an ideal $I$. To make sure the ideal $I$ is minimal, we have to choose carefully the positions of $\{\l_i-1, \dots, 1-\l_i\}$ in the sequence  $(h_1, \dots, h_n)$.

Similar ideas apply to  other types. For classical groups  of types $B, C $ and $D$, the construction of  the weighted Dynkin diagram is a little different and we have to adjust our choice accordingly.


\section{Minimal Ideals For Type $A_{n-1}$}

Suppose $\ccg= \ccs \ccl(n)$. Following standard notation, let
$\ccb$   be  the standard upper triangular matrices and
$\cch$   be the diagonal matrices. The root system  of $\ccg$ is
$\{e_i- e_j \mid 1 \le i,j \le n, i\neq j \} $ and $\D^+=\{e_i- e_j
\mid 1\le i <j \le n\}$. We denote by $E_{ij} $ the elementary matrix with its $ij$-entry 1 and other entries 0.  The  root space for $e_i- e_j$ is spanned by the matrix $E_{ij}$.

By Theorem \ref{cm1}, let $\l=[\l_1,\l_2, \dots, \l_p]$ be a partition of $n$. Following  the  idea discussed at the end of  previous section, we need to construct some maps to keep track of the positions of  $\{\l_i-1, \l_i-3, \dots, 1-\l_i\}$.


Indeed, let $\s_{i}$, for $1 \le i \le p$,  be a sequence of index maps $$ \s_{i}: \{\l_i-1, \l_i -3, \dots, 1-\l_i\} \rightarrow  [n]=\{1,2, \dots, n\}$$

\begin{lem}\label{sigma} There exists a sequence of maps $ \{ \s_{i} \}_{i=1}^p$, satisfying the following  properties:



(1)  Each $\s_{i}$ is one-to-one and $Im(\s_i) \cap Im(\s_j)= \emptyset$, if $i \neq j$.

(2)  If $k < l$ and $k \in Dom(\s_i)$, $l \in Dom(\s_j)$, then $\s_{i}(k) > \s_{j}(l)$.

(3)  For any $\l_i, \l_j$ and $k, l \in Dom(\s_i) \cap Dom(\s_j)$, if $\s_{i}(k) > \s_{j}(k)$, then $\s_{i}(l) > \s_{j}(l)$.
Here  $Dom(\s_{i})$ denotes the domain of $\s_{i}$ and $Im(\s_{i})$  denotes the image of $\s_{i}$.
\end{lem}

Proof. We form a sequence of integers $h=(h_1, \dots, h_n)$ by placing $\l_i-2s+1$ in the position $\s_i(\l_i-2s+1)$ for $1\le i \le n$. Property (1) and (2) make sure  we indeed get a weighted Dynkin diagram from $h$. Property $(3)$ gives some restriction on the positions of integers of the same value, but coming from different parts of $\l$. \qed

\begin{rmk}
1. If $\{ \s_1, \dots, \s_p\}$ is a sequence of maps as above, then
$\sqcup_{i=1}^p Im(\s_{i})=[n]$.

2. As a special case of property (2), $\s_{i}(\l_i -1) < \s_{i}(\l_i -3) < \dots < \s_{i}(1-\l_i )$.
\end{rmk}

\begin{example}\label{eg} Let $\l=[4, 2]$ be a partition of $6$. Let $h= (\underline 3, \underline 1, \overline 1, \underline {-1}, \overline {-1}, \underline {-3})$ and $h'= (\underline 3, \overline 1, \underline  1, \overline {-1}, \underline {-1}, \underline {-3})$, where $\overline i$ means that $i$ comes from  $\l_2=2$ and $\underline i$ means that $i$ come from $\l_1=4$ of the partition $\l$. Then  $h, h'$ give rise to the same weighted Dynkin diagram and also show that the sequence of the maps $\{\s_i \}$ is not unique.
\end{example}

For each $\s_i$, we attach a set of positive roots:  $$\cc^+ (\s_i)= \{ e_{\s_{i}(\l_i -1)}- e_{\s_{i}(\l_i -3)},  e_{\s_{i}(\l_i -3)}- e_{\s_{i}(\l_i -5)}, \dots, e_{\s_{i}(3-\l_i)} -e_{\s_{i}(1-\l_i)}\}.$$

For the partition $\l=[\l_1, \l_2, \dots, \l_p]$, we set $\cc= \cup_i \cc^+(\s_i)$. Let $X_\a$ be a root vector that corresponds to the root $\a \in \cc$ and define $X$ be the sum of $X_\a$.

Given $\l$, attach $\cc_+(\l_i)=\{e_{N_i+1}- e_{N_i+2}, \dots, e_{N_i+\l_i-1}-e_{N_i+\l_i}\}$ to $\l_i$. Here $N_i$ is chosen so that $\sqcup_{ 1 \le i \le p, 1\le j\le \l_i}(N_i+j)=[n]$. Let $\cc_+$ be the union of $\cc_+(\l_i)$ and let $\tilde X= \sum_{\a \in \cc_+}X_\a$. It's showed in  \cite[section 5.2]{CM} that $\tilde X$ lies in the orbit $\co_\l$. From the construction of $\cc$ and $\cc_+$, we can see that $X$ is conjugate to $\tilde X$ by some elements in  $S_n$, therefore also lies in $\co_\l$.

Now let $$H_{\cc^+ (\s_i)}= \sum_{s=1}^{\l_i}(\l_i -2s +1 )E_{\s_{i}(\l_i-2s +1), \s_{i}(\l_i-2s +1)} $$ and $$ H= \sum _{i=1}^{n} H_{\cc^+(\s_i)}.$$

Although each $H_{\cc^+(\s_i)}$ is dependent on the choice of the map $\s_{i}$, by property (2) above, the diagonal entries of $H$ are decreasing and $H$ is independent of the series  of the maps $\{\s_{i}\}$ we choose. Indeed, $H$ is  the Dynkin element of the orbit $\co_\l$.

For any root $\a=e_{\s_{i}(k)}-e_{\s_{j}(l)} \in \D$, $\a(H)=k-l$. In particular \[\tag{2.2.1} \a(H)=2, \qquad \text{ for any } \a \in \cc.\]

There are several things to show:

\begin{lem}\label{chain}
The set of positive roots $\cc$ is an antichain in $\D^+$.
\end{lem}

Proof. The set  $\cc$ is an antichain if and only if
 $\a-\b \notin Q^+$ for any roots $\a, \b \in \cc$.  Let $\a= e_{\s_i}(m)-
e_{\s_i}(m-2)$ and $\b= e_{\s_j}(l)- e_{\s_j}(l-2)$ be two roots in $\cc$. Then by
property (3) in lemma \ref{sigma} above,  $\a-\b= (e_{\s_i}(m)-  e_{\s_j}(l))+
(e_{\s_j}(l-2)- e_{\s_i}(m-2))$ can not lie in $Q^+$. Thus $\cc$ is an antichain. \qed

Then the set $\cc$ satisfies the assumption of Lemma \ref{trip} and we can find an appropriate nilnegative element $Y$ such that $\{H, X, Y\}$ is a standard triple.

Since there is a canonical bijection between the  antichains of the root poset and the ad-nilpotent ideals,  we  can construct an ad-nilpotent ideal $I_\cc$ that is generated by  $\cc$.

\begin{lem}\label{h3}
The ideal $I_{\cc}$ contains $\ccq_{H,3}$.
\end{lem}

Proof. Both $I_\cc$ and $\ccq_{H, 3}$ are direct sum of root spaces. Let $\a$ be a positive root and $\ccg_\a \subset \ccq_{H,3}$.
Suppose $\a= e_{\s_{i}(k)}-e_{\s_{j}(l)}$. Then $\a(H)= k-l \ge 3$.
There are two possible cases for $k$.  If $k >0$, then $k- 2 \in
Dom(\s_{i})$ and $\a =\b +\g$, where $\b=
e_{\s_{i}(k)}-e_{\s_{i}(k-2)} $ and $\g= e_{\s_{i}(k-2)}-
e_{\s_{i}(l)} $.  By  property 2  of the maps $\{\s_i$, $\s_j \}$,
$\b \in  \cc$ and $\g \in \D^+$, hence $\a \in \D(I_{\cc})$. If $k
\le 0$, then $l<0$ and $l+2 \in Dom(\s_{j})$.  Therefore
$\a= \b +\g \in \D(I_\cc)$, where $\b=e_{\s_{i}(k)}- e_{\s_{i}(l+2)}
\in \D^+ $ and $\g= e_{\s_{i}(l+2)}-e_{\s_{i}(l)} \in \cc$.  \qed

\begin{prop}\label{min}
The associated orbit of the ideal $I_\cc$ is  $\co_\l$ and it is an ideal of minimal  dimension.
\end{prop}

Proof. By formula $(2.2.1)$, the ideal $I_\cc$ is contained in the Dynkin ideal $\ccq_{H,2}$. By Kostant's theorem \ref{Kos}, the associated orbit of the ideal $I_\cc$ is contained in the closure of the orbit $\co_X $. On the other hand, $I_\cc$ contains $X$, so $\co_{I_\cc}=\co_X$. We only need to prove that $\dim I_\cc=m_\co$.

The Dynkin ideal $\ccq_{H,2}$ is contained in the Borel subalgebra $\ccb$ and there is a  decomposition $\ccb= \ccg^+_0 \oplus \cch \oplus \ccg_{H,1} \oplus \ccq_{H,2}$. By  the formula  2.1.1 for $m_\co$, \begin{align*} \dim \ccq_{H,2}-m_{\co} &=\dim B - (\dim \ccg_{H,1} +\dim \ccg^+_0 + \dim \cch)-\dim B+\dim B_{G_X} \\ &=\dim B_{G_X}- (\dim \ccg_{H,1} +\dim \ccg^+_0 + \dim \cch) \\ &= \frac{1}{2}[\dim (\ccg_{H, 0})+ \dim(\ccg_{H, 2})+\text{ rank }G_X]-\dim \ccg^+_0- \dim \cch \\ &=\frac{1}{2} (\dim \ccg_{H,2}- \dim \cch + \text{rank } G_X ).\end{align*}  The last equality follows from the fact that $\ccg_{H, 0}= \ccg_0^+ \oplus \ccg_{0}^- \oplus \cch$ and $\dim \ccg_0^+ =\dim \ccg_0^-$.


Let $\cc^-=\{ \a \in \D(\ccg_{H,2}) \mid \a \notin \D(I_\cc) \}$ and $\cc^+=\{ \a \in  \D(\ccg_{H,2}) \cap \D(I_\cc) \mid \a \notin \cc \}$.

Then $\D(\ccg_{H,2})=  \cc \sqcup \cc^+ \sqcup \cc^- $.    By Lemma \ref{pa}, since the set $\cc$ consists of linearly independent roots, $|\cc|= \dim ad_X(\cch)= \dim \cch - \dim Z_{\cch}(X)= \dim \cch- \text{rank }G_X $. Then $$\dim \ccq_{H,2}- m_\co=\frac{1}{2}(|\cc^+| + |\cc^-|).$$

From Lemma \ref{h3}, $\ccq_{H,3}$ is contained both in $\ccq_{H,2}$ and $I_\cc$.  So \[\dim \ccq_{H,2}- \dim I_\cc= \dim \ccg_{H,2}-\dim I_\cc \cap \ccg_{H,2}=|\cc^-|. \]

Now it suffices to  prove that $\cc^+$ and $\cc^-$ have the same cardinality. This follows from Lemma \ref{inv} below, which proves this proposition. \qed

\

Notice that any root in  $\D(\ccg_{H,2})$ has the form  $e_{\s_{i}(m)}-e_{\s_{j}(m-2)}$ for some $i, j$ and $m \in Dom(\s_{i})$ and $ m-2 \in Dom(\s_{j})$. We define a map $$\iota: \D(\ccg_{H,2}) \rightarrow \D(\ccg_{H, 2})$$ by $e_{\s_{i}(m)}-e_{\s_{j}(m-2)} \mapsto e_{\s_{j}(2-m)}-e_{\s_{i}(-m)}$.  Since the domain of $\s_{i}$ and $\s_{j}$ is symmetric with respect to 0, so $-m\in Dom(\s_{i})$ and $2-m \in Dom(\s_{j})$ and the map is well-defined.

\begin{lem}\label{inv} Keep the notations as above, the map $\iota$ is an involution on $\D(\ccg_{H,2})$. Moreover, $\iota$ maps $\cc$ to itself and maps $\cc^+$ to $\cc^-$ and vice versa.
\end{lem}

Proof. Since $\iota^2=id$, it's obvious that $\iota$ defines an involution on $\D(\ccg_{H,2})$. If $\a=e_{\s_{i}(m)}-e_{\s_{j}(m-2)} \in \cc$, then $m, m-2$ come from the same part of the partition $\l$, so $i =j$ and $\iota(\a) \in \cc$.  That means $\iota$ maps $\cc$ to itself and also maps $\cc^+ \sqcup \cc^-$ to itself.

Suppose that $\a \in \cc^+$. Since $\a$ lies in the ideal that's generated by $\cc$, there exists a root $\b \in \cc$ and a positive root $\g \in \D^+$ such that  $\a =\b+ \g$.

If $m >0$, then $m-2 \in Dom(\s_{i})$. In this case $\a= \b +\g$, where $\b= e_{\s_{i}(m)}- e_{\s_{i}(m-2)} \in \cc  $ and $\g= e_{\s_{i}(m-2)}- e_{\s_{j}(m-2)}>0 $.  Hence $\s_{i}(m-2)< \s_{j}(m-2) $.  The domains of $\s_{i}$ and $\s_{j}$ are symmetric with respect to 0, so $2-m, -m \in Dom(\s_{i})$ and $2-m \in Dom(\s_{j})$. Then   $\iota(\a)= e_{\s_{j}(2-m)}-e_{\s_{i}(-m)} = \b' + \g'  $ , where $\b'= e_{\s_{j}(2-m)}-e_{\s_{i}(2-m)} $ and $\g'= e_{\s_{i}(2-m)}-e_{\s_{i}(-m)}  $. Then $\g' \in \cc$ and by property (3) of the maps $\{\s_i, \s_j\}$, $\s_{i}(2-m) < \s_{j}(2-m)$, so $\b' \in \D^- $ and $\iota(a) \in \cc^-$.

If $m \le 0$, then $m \in Dom(\s_{j})$. In this case, $\b= e_{\s_{j}(m)}-e_{\s_{j}(m-2)} \in \cc  $. With the same argument, $\a'= e_{\s_{j}(2-m)}-e_{\s_{i}(-m)}$ is the unique element that corresponds to $\a$ and lies in $  \cc^-    $.

By the same reasoning, $\iota$ maps $\cc^-$ to $\cc^+$. The lemma is proved.
\qed

\

\section{Dimension formula for minimal ideals of Type $A_{n-1}$}

For type $A_{n-1}$, $\ccn$ is the set of strictly upper triangular
matrices. Following  \cite{Pa1}, an ad-nilpotent ideal is
represented by  a right-justified Ferrers (or Young) diagram with
at most $n-1$ rows, where the length of the $i$-th row is at most
$n-i$. Namely, if any root space $\ccg_\a$ lies in an  ideal $I$,
then any root subspace $\ccg_{\b}$ that's on the northeast side of
$\ccg_{\a}$ also lies  this ideal. The generators of the ideal are
the set of   southwest corners of the diagram. We use the pair $[i,
j]$ to denote the positive root $e_i-e_j$, where $1 \le i <j \le n
$. Then the set of generators of the ideal $I$ can be writen as (see Figure 2-1):
\begin{align*}\G(I)=\{[i_1, j_1], \dots, [i_k, j_k]\},  \text{where
} & 1 \le i_1 <\dots <i_k \le n-1, \\ & 2 \le j_1 <\dots<j_k \le n.
 \end{align*}

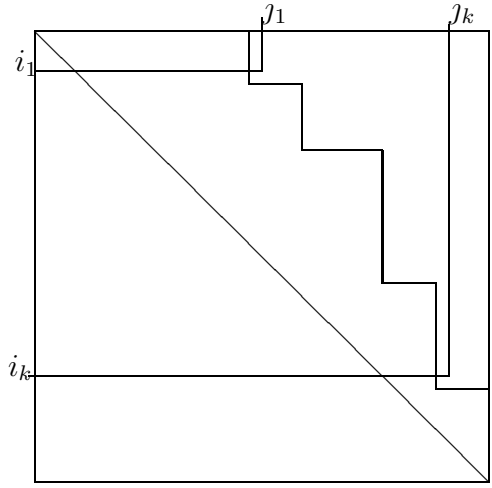
\begin{figure}
\setlength{\unitlength}{5pt}
\begin{center}
\begin{picture}(40, 40)(0, 0)
\put(2, 2){\line(1, 0){34}}
\put(2, 2){\line(0, 1){34}}
\put(2, 36){\line(1, 0){34}}
\put(36, 2){\line(0, 1){34}}
\put(2, 36){\line(1, -1){34}}
\put(18, 36){\line(0, -1){4}}
\put(18, 32){\line(1, 0){4}}
\put(22, 32){\line(0, -1){5}}
\put(22, 27){\line(1, 0){6}}
\put(28, 27){\line(0, -1){10}}
\put(28, 17){\line(1, 0){4}}
\put(32, 17){\line(0, -1){8}}
\put(32, 9){\line(1, 0){4}}
\put(19, 33){\line(-1, 0){17}}
\put(19, 33){\line(0, 1){4}}
\put(0.5, 33){$i_1$}
\put(19, 37){$j_1$}

\put(33, 10){\line(-1, 0){31.5}}
\put(33, 10){\line(0, 1){26.5}}
\put(33, 37){$j_k$}
\put(0, 10){$i_k$}

\end{picture}

\end{center}
\vspace{-.5cm}
\begin{quote}
\caption{An ad-nilpotent ideal for type $A_{n-1}$
\label{fig1}}
\end{quote}
\vspace{-1cm}
\end{figure}

In order to compare the dimension of two minimal ideals, we first need to have an explicit formula for the dimension of the minimal ideal in terms of  the partition $\l$. Suppose the sequence of maps $\{\s_{i}\}$ satisfies   properties (1) and (2) in Lemma \ref{sigma}  plus an additional one:

(4)  $\s_{i}(k) <\s_{j}(k)$, when  $i<j$  and  $k$  lies in the domain of  $\s_{i}$  and  $\s_{j}$.

The sequence of the maps  exists and  is uniquely determined by those restrictions. Moreover,  $\{\s_i \}$ automatically satisfy property (3) of Lemma \ref{sigma} so the ideal $I$
constructed  from these maps has  minimal dimension.   The
Dynkin element  $H= diag\{h_1, h_2,\dots  h_n\}$ is the same as in last section.

Let $A(l)$ be the number of entries in $H$ that are less or equal to $l$. Let $B(l)$ be the number of entries of $H$ that are bigger than $l$. Then $A(l)+B(l)=n$.  If we rewrite the partition $\l=[\l_1, \l_2, \dots , \l_p]$ in the exponential form $\mathbf{\l}= [t_1^{n_{t_1}}, t_2^{n_{t_2}}, \dots, t_k^{n_{t_k}}]$, then the minimal dimension has the following formula: \begin{prop}
The  dimension of the minimal ideals corresponding to  the partition $\l$  is equal to
$$m_{\co_\l}= \frac{n(n+1)}{2} +\sum_{i=1}^n (-A(\l_i-1)+A(-\l_i-1))+\sum_{i=1}^k \frac{n_{t_i}(n_{t_i}-1)}{2}.$$
\end{prop}

Proof: To calculate the dimension for the ideal $I$, we need to sum up the number of positive  roots in  $I$ in  each  row.

From the construction of the maps $\{\s_{i}\}$, it is obvious that  in the $\s_{i}(\l_i-1), \s_{i}(\l_i-3) \dots \s_{i}(3-\l_i)'s$ rows, $\cc_{\l_i}$ forms a subset of generators of the ideal $I$, therefore the Ferrers diagram begins with boxes ${[\s_{i}(\l_i-1), \s_{i}(\l_i-3)], \dots, [\s_{i}(3-\l_i), \s_{i}(1-\l_i)]} $.

On the other hand, there is no generator in row $\s_{i}(1-\l_i)$.  Let $[s_i, t_i]$ be the generator that's below this row and row $\s_{i}(1-\l_i)$.  Then row $s_{\l_i}$ share the same columns. Then row $\s_i(1-\l_i)$ of the Ferrers diagram begins with box $[\s_{i}(1-\l_i), t_i]$. If there is no generator below this row, we simply say that the diagram begins with box $[\s_{i}(1-\l_i), n+1]$. This convention makes the formula $(2.3.1)$ in the next paragraph give the correct number of positive roots in row $\s_i(1-\l_i)$, which is zero. 

Suppose that the Ferrers diagram  corresponding to the ideal $I$ begins with box $[s, t]$ in row $s$. Then the number of positive roots in this row is equal to $1+n-t=(1+n-s)+(s-t)$.

The total summation of positive roots in all  rows is equal to:  \begin{align*}\tag{2.3.1}
|I| &= \sum_{i=1}^p \sum_{s=2}^{\l_i}\big[ n+1-\s_{i}(\l_i-2s+1)\big]+ \sum_{i=1}^p \big[ n+1-t_i \big]
\\ &=\sum_{i=1}^p \sum_{s=1}^{\l_i}(n+1-\s_{i}(\l_i-2s+1))+ \sum_{i=1}^p\s_{i}(\l_i-1)-\sum_{i=1}^p t_i.
\end{align*}

Since $\cup_{i=1}^p Im(\s_{i})=[n]$, the first term in the equation  $(2.3.1)$ is equal to $\sum_{j=1}^n(j)=\frac{n(n+1)}{2}$. We only need to know the value of $\s_{i}(\l_i-1)$ and $t_i$.

Indeed, $\s_{i}(\l_i-1)$ shows the position of $\l_i-1$ in the Dynkin element $H=diag\{h_1, h_2,\dots  h_n\}$. If $\l_i$ satisfies that all $\l_j<\l_i$, when $j>i$, by property 4, $\s_{i}(\l_i-1)$ is the  first $t$, such that $h_t=\l_i-1$. Namely $\s_{i}(\l_i-1)= max\{t \mid  h_t  >\l_i-1\}+1= B(\l_i-1)+1$.

If $\l_{i-1}=\l_i$, then $\s_{{i-1}}(\l_{i-1}-1)=\s_{i}(\l_i-1)+1$. If we consider the exponential expression of the partition $\l$, the summation of all $\s_{{i}}(\l_i-1)$, where $\l_i=t_j$ is equal to $n_j(1+B(\l_i-1))-n_j(n_j+1)/2$.

What remains to discuss is the value of $t_i$. We need to find the  nearest corner of the Ferrers diagram that's below row $\s_{{i}}(1-\l_i)$.

Case 1: Suppose that $\l_l$ is the smallest integer such that $\l_l> \l_i$ and $\l_l \equiv \l_i (\text{mod } 2)$. Then $[\s_{{l}}(1-\l_i), \s_{{l}}(-\l_i-1)] $ is a generator that's below row $\s_{{i}}(1-\l_i)$. The column coordinate $t_i$ is the smallest integer $\s_{{l}}(\l_i-1)$ for such $\l_l$. Therefore it is equal to $B(-\l_i-1)$.

Case 2: Suppose there's no such $\l_l$ as in case (1). But there exist some $\l_l$ such that    $\l_l \ge  \l_i+2 $.  In this case $[\s_{{i}}(-\l_i), \s_{{i}}(-2-\l_i)]$ is a generator that's below row  $\s_{{i}}(1-\l_i) $. Then similar to case (1), $t_i$ is equal to $B(-2-\l_i)+1$. But  under previous assumption, no $-1-\l_i$ appears in the diagonal entries of $H$, so $B(-2-\l_i)=B(-1-\l_i)$.

Case 3: Suppose $\l_l -\l_i \le 1$ when $1 \le l \le i$. Then all diagonal entries of $H$ are bigger than $-l_i-1$ and  $B(-1-\l_i ) =n$. In this case, there's no generator below row $\s_{{i}}(1-\l_i)$ so $(\star)_{\l_i}= n+1=B(-1-\l_i+1)$.

The formula for $m_{\co_\l}$ is derived if we  use $A(l)=n- B(l)$ and put the values of $\s_{{i}}(1-\l_i)$ and $t_i$ into the equation $(2.3.1)$. \qed

\

\begin{lem}\label{al}
$A(l)=\sum_{i=1}^n max(min(\lfloor \frac{\l_i+l+1}{2} \rfloor, \l_i),0)$.
\end{lem}

Proof. The number of elements in the set $\{\l_i-1, \dots, 1-\l_i \}$ that are at most $l$ is equal to a positive integer $t$, where $t \le \l_i$ and $-\l_i-1 +2t \le l$. The summation of all such $t$ is $A(l)$. \qed

\

We write $\co_1 \le \co_2$ (resp. $\co_1 < \co_2 $) if the closure of the orbit $\co_1$ is (resp. strictly) contained in the closure of the  orbit $\co_2$. This defines a partial  order on nilpotent orbits. It is obvious that if $\co_1$ is smaller than $\co_2$, the dimension of $\co_1$ is smaller than the dimension of $\co_2$. It turns out that we also have  the same relation for the dimension  of   minimal ideals.

Suppose that ${\l} =[\l_1, \l_2, \dots \l_p]$ and  $\bold d =[d_1, d_2, \dots, d_q]$ are two partitions of $n$ and correspond to the orbits $\co_{\l}$ and $\co_{\bold d}$ respectively. As shown in  \cite[6.2.1]{CM}, the partial order on $P(n)$ is defined as:  $\bold d \le \l$  if and only if $\sum_{1 \le j \le l}\l_j \le \sum_{1 \le j \le l} d_j$ for $1 \le l \le n$.

\begin{lem}  \cite[Lem6.2.4]{CM}\label{cm} (1). Suppose $\l, \bold d \in P(n)$. Then  $\l$ covers $\bold d$ in the order $\le$ (meaning $\l<d $ and there is no partition $e$ with $\l < e <d$) if and only if $\bold d$ can be obtained from $\l$ by the following procedure. Choose an index $i$ and let $j$ be the smallest index greater than $i$ with $0 \le \l_j < \l_i-1$. Assume that either $\l_j=\l_i-2$ or $\l_k=\l_i$ whenever $i <k< j$. Then the parts of $\bold d$ are obtained from the $\l_k$ by replacing $\l_i, \l_j$ by $\l_i-1$, $\l_j+1$ respectively(and rearranging).

(2). $\co_\l \le \co_{\bold d} $ if and only if $\l \le \bold  d$. Hence $\l $ covers $\bold d$ iff $\co_{\bold d} <\co_{\l}$ and there is no nilpotent orbit $\co_e$, with $\co_{\bold d} <\co_e <\co_\l$.

\end{lem}

 \begin{prop}\label{strict}
If $\co_{\bold d} \le \co_{\l}$ (resp $\co_{\bold d} < \co_\l $), then $m_{\co_{\bold d}} \le m_{\co_\l}$  (resp. $m_{\co_{\bold d}} < m_{\co_\l}$).
 \end{prop}

Proof. It suffices to prove the proposition under the assumption
that $\l$ covers $\bold d$. Since all $\l_i$ are nonnegative, by Lemma \ref{al},
$$A(\l_i-1)=\sum_{j=1}^n min(\lfloor \frac{\l_i+\l_j}{2} \rfloor,
\l_j)=\sum_{j \le i} \lfloor \frac{\l_i+\l_j}{2} \rfloor+ \sum_{j <
i}\l_j$$
$$A(-1-\l_i)=\sum_{j}max(min( \lfloor \frac{\l_j-\l_i}{2} \rfloor,j), 0)=\sum_{j<i} \lfloor \frac{\l_j-\l_i}{2} \rfloor.$$
Therefore
\begin{align*}
A(\l_i-1)-A(-\l_i-1) &= \sum_{j<i}(\lfloor \frac {\l_i+\l_j}{2} \rfloor-\lfloor \frac{\l_j-\l_i}{2} \rfloor)+ \sum_{j \ge i}\l_j
\\ &=\sum_{j<i}\l_i +\sum_{j \ge i} \l_j.
\end{align*}

Hence,

$$\sum_{i}(\sum_{j < i}\l_i+ \sum_{j \ge i}\l_j)=\sum_{i}(\sum_{j <i}1)\l_i+ \sum_{j}(\sum_{i \le j }1) \l_j =\sum_{i}(2i-1)\l_i.$$

Therefore the formula for $m_\co$ can be written as $$m_\co = \frac{n(n+1)}{2} +\sum_{i=1}^n(2i-1)\l_i + \sum_{i=1}^k \frac{n_{t_i}n_{t_i}-1}{2}.$$

Since $\l$ covers $\bold d$, as in lemma \ref{cm},  for simplicity, we assume when $i<t<j$, $\l_i < \l_t<\l_j$.  Then $d$ differs from  $\l$ only when $d_i = \l_i-1$ and $d_j= \l_j+1$. The difference between the second term of  $m_{\co_\l}$ and $m_{\co_{\bold d}}$ is equal to
\begin{align*}
& -(2i-1)\l_i-(2j-1)\l_j+(2i-1)(\l_i-1)-(2j-1)(\l_j+1) \\ &=2(j-i) >0.
\end{align*}

We need to discuss the exponential  form of $\l$ and $\bold d$.  If
$\l_i-\l_j >2$, then because of the assumption, $i=j-1$. Since $\l_i
\neq \l_j$, suppose the exponential expression of $\l$ is
$\l=[\l_1^{n_{\l_1}}, \dots, \l_i^{n_{\l_i}}, \l_j^{n_{\l_j}},
\dots, \l_{n}^{n_{\l_n}}]$. Then the exponential form of $\bold d$ is
$\l=[\l_1^{n_{\l_1}}, \dots, \l_i^{n_{\l_i}-1}, \l_i-1, \l_j+1,
\l_j^{n_{\l_j}-1}, \dots, \l_{n}^{n_{\l_n}}]$. The difference
between the third term of $m_{\co_\l}$ and $m_{\co_{\bold d}}$ is
\begin{align*}
&\frac{n_{\l_i}(n_{\l_i}-1)}{2}+\frac{n_{\l_j}(n_{\l_j}-1)}{2}-\frac{(n_{\l_i}-1)(n_{\l_i}-2)}{2}-\frac{(n_{\l_j}-1)(n_{\l_j}-2)}{2}
\\ &= n_{\l_i}+n_{\l_j}-2 \ge 0.
\end{align*}

From the two inequalities above, it's easy to deduce that $m_{\co_\l} > m_{\co_{\bold d}}$.

If $\l_i-\l_j =2$, then $\l_t =\l_i-1 = \l_j+1$, for any $i<t<j$. The exponential form of $\l$ is $\l=[\l_1^{n_{\l_1}}, \dots, \l_i^{n_{\l_i}}, \l_t^{n_{\l_t}},  \l_j^{n_{\l_j}}, \dots, \l_{n}^{n_{\l_n}}]$ and $d$ has exponential form $d=[\l_1^{n_{\l_1}}, \dots, \l_i^{n_{\l_i}-1},  \l_t^{n_{\l_t}+2}, \l_j^{n_{\l_j}-1}, \dots, \l_{n}^{n_{\l_n}}]$. The difference between the third term of $m_{\co_\l}$ and $m_{\co_{\bold d}}$ is

\begin{align*}
&\frac{n_{\l_i}(n_{\l_i}-1)}{2}+\frac{n_{\l_j}(n_{\l_j}-1)}{2}+\frac{n_{\l_t}(n_{\l_t}-1)}{2} \\& -\frac{(n_{\l_i}-1)(n_{\l_i}-2)}{2}-\frac{(n_{\l_j}-1)(n_{\l_j}-2)}{2}-\frac{(n_{\l_t}+2)(n_{\l_t}+1)}{2} \\ &=n_{\l_i}+n_{\l_j}-2 -2n_{\l_t}-1. \end{align*}

Since $n_{\l_t}= j-i-1$, we can compare the second and  third term of $m_{\co_\l}$ and $m_{\co_{\bold d}}$ and still get strict inequality. \qed

\

We proved the existence of  minimal ideals of dimension $m_\co$ for nilpotent $\co$. And Example \ref{eg} shows that minimal ideals are not  unique. All the minimal ideals  we have constructed above are contained in the Dynkin ideal $\ccq_{H,2}$.  However, it's possible to have   minimal ideals that are not contained in the Dynkin ideal. If we know some information about a general  ideal $I$, here is a criterion to see whether this ideal is minimal or  not.

\begin{cor}
Suppose the ideal $I$ contains a nilpotent element $X$ and $dim I=m_{\co_X}$, then the associated orbit of $I$ is $\co_X$.
\end{cor}

Proof. The ideal $I$ contains $X$, so  $\co_I \ge \co_X$. By the strict inequality in Proposition \ref{strict},  it's not possible that  $\co_I \gneq \co_X$. \qed

\

\section{Minimal Ideals For Type $C_n$}

Let $\ccg$ be  $\ccs \ccp_{2n}$. The set of positive roots in $\ccg$ is $\D^+=\{e_i-e_j \mid 1 \le i<j \le n \} \cup \{2e_i \mid 1 \le i \le  n \}$. Take a partition $ \bold d$ of $2n$ with odd parts repeated with even multiplicity. By Theorem \ref{cm1}, it corresponds to a nilpotent orbit $\co_{\bold d}$ in $\ccg$. The partition $\bold d$ can be written as $$\bold d=[d^{r_d}, (d-1)^{r_{d-1}} , \dots, 2^{r_2}, 1^{r_1}],$$  where  $r_k$  is even  if $k$ is odd.  To distinguish the integers with  with the same value but in different positions,  we rewrite $\bold d$ as $$\bold d=[d_1, \dots, d_{r_d}, \dots, 1_{1}\dots , 1_{r_1} ],$$ where  $ \sum_{i=1}^i (r_i) d=2n$  and  $r_i \text{ is even if } i \in 2\mathbb N+1$.  Here the index $d_i$ has value $d$ and the subscript $i$ distinguishes indexes with the same value but at different positions.

The procedure to  get an ideal of minimal dimension for the nilpotent orbit $\co_{\bold d}$ is similar to  what  we did in the type $A_{n-1}$ case.  Let $\ca$ be the index set $\sqcup_{k=1}^{d} \{ k_1, \dots, k_{r_k} \}$. First recall the procedures to get the weighted Dynkin diagram of $\bold d$. We take the union of integers $\{k_i-1, k_i-3,  \dots, 1-k_i\}$ for any $k_i \in \ca$ and rearrange the sequence in the form $h= (h_1, \dots, h_n, -h_1, \dots, -h_n)$, where $h_1 \ge h_2 \dots \ge h_n \ge 0$.  Again we need some index maps to keep track of the first $n$ integers in $h$. 

We define the following  sequence of  maps  $\{\s_{k_i}\}_{k_i \in \ca}$. 

\[
\s_{k_i}=\begin{cases}
   \{ k-1, k-3, \dots, 2, 0\} \rightarrow [n] \text{ if } k \text{ is odd and } 1 \le i \le \frac{r_{k}}{2}; \\
    \{ k-1, k-3, \dots, 2\} \rightarrow [n] ,\text{ if } k \text{ is odd and }\frac{r_{k}}{2} < i \le r_{k}; \\
    \{ k-1, k-3, \dots, 1\} \rightarrow [n], \text{ if } k \text{ is even and }1 \le i \le r_{k}.
    \end{cases}
\]
For any $k_i \in \ca$, let  $Im(\s_{k_i})$ (resp. $Dom(\s_{k_i})$) be the image (resp. domain) of $\s_{k_i}$. For $1 \le k \le 2n, 1 \le i \le  r_k$, write $\tilde i_k = r_k+1 -i$. 

\begin{lem}\label{cnmap}
There exists a set of maps $\{ \s_\t\}_{\t \in \ca}$ satisfying the following properties:

(1) For any $\t, \o \in \ca$, $\s_{\t}$ is one-to-one and $Im(\s_\t) \cap Im(\s_\o) =\emptyset$.

(2) For any $\t, \o \in \ca$, $m \in Dom(\s_\t)$, $l \in Dom(\s_\o)$, and $m >l$, $\s_\t(m) < \s_\o(l).$

(3) If $1 \le k \le 2n$ , $1 \le i< j \le r_{k}$, and $m \in Dom(\s_{k_i}) \cap  Dom(\s_{k_j})$, then $\s_{k_i}(m)< \s_{k_j}(m)$.

(4) If $k, l \in \mathbb N$ , $1 \le i \le \lfloor \frac{r_k}{2} \rfloor$, $ 1\le j \le \lfloor \frac{r_l}{2} \rfloor$ and $m > 0$, then either  $\s_{l_{j}}(m)<\s_{k_{i}}(m)< \s_{ k_{\tilde i_{k}}}(m)< \s_{l_{\tilde j_l}}(m)$ or $\s_{k_{i}}(m)<\s_{l_{j}}(m)< \s_{ l_{\tilde j_l}}(m)< \s_{ k_{\tilde i_k}}(m)$.

(5) Let $k, l,  i, j$ be the same as in (4), if $\s_{k_i}(2)< \s_{l_j}(2)$, then  $\s_{k_i}(0)< \s_{l_j}(0)$.

(6) Let $k, l$ be even integers, if $r_k$ is odd, $i=\lceil \frac{r_k}{2} \rceil $ and $ 1\le  j \le \lfloor \frac{r_l}{2}\rfloor$, then $\s_{l_j}(m)< \s_{k_i}(m)< \s_{ l_{\tilde j_l}}(m)$.


(7) Let $k,l$ be even integers. If $r_k$, $r_l$ are odd and $k<l$, then $\s_{k_{\lceil \frac{r_k}{2} \rceil}}(m)<\s_{l_{\lceil \frac{r_l}{2} \rceil}}(m)$.

\end{lem}

The construction of index maps is similar to type $A_{n-1}$ case except that here we need more restrictions for different types of indexes maps. If we put any $m \in Dom(\s_\t)$ into the position $\s_\t(m)$ and $-m$ into $\s_\t(m)+n$, properties (1) and (2) make sure that we could get $h$ as above. Property (3) and (4) gives orders for integers from different parts of $\bold d$. Properties  5 deal with odd parts of the partition $\bold d$. And properties (6) and (7) deal with even parts of $\bold d$.


Now it's possible to get a set of positive roots. If $k$ is odd and $1 \le i \le \frac{r_k}{2}$, set $$\cc(\s_{k_i})= \{ e_{\s_{k_i}(k-1)}-e_{\s_{k_i}(k-3) }, \dots, e_{\s_{k_i}(2)}- e_{\s_{k_i}(0)}\}. $$

If $k$ is odd and $i>\frac{r_k}{2}$, set $$\cc(\s_{k_i})= \{ e_{\s_{k_i}(k-1)}-e_{\s_{k_i}(k-3) }, \dots, e_{\s_{k_i}(4)}- e_{\s_{k_i}(2)}, e_{\s_{k_i}(2)}+e_{\s_{ k_{\tilde i_k}}(0)}\}. $$

If $k$ is even and $i \le \lfloor \frac{r_k}{2} \rfloor$, set $$ \cc(\s_{k_i}) =
\{ e_{\s_{k_i}(k-1)}- e_{\s_{k_i}(k-3)} , \dots,  e_{\s_{k_i}(3)}- e_{\s_{k_i}(1)}, e_{\s_{k_i}(1)}+ e_{\s_{ k_{\tilde i_k}}(1)} \}. $$

If $k$ is even and $i>\lceil \frac{r_k}{2} \rceil$, set $$ \cc(\s_{k_i}) =
\{ e_{\s_{k_i}(k-1)}- e_{\s_{k_i}(k-3)} , \dots,  e_{\s_{k_i}(3)}- e_{\s_{k_i}(1)} \}. $$

If $k$ is even, $r_{k}$ is odd and $i=\lceil  r_k/2 \rceil$, set
$$\cc(\s_{k_i})=\{e_{\s_{k_i}(k-1)}-e_{\s_{k_i}(k-3)}, \dots, e_{\s_{k_i}(3)}-e_{\s_{k_i}(1)}, 2e_{\s_{k_i}(1)}\}.$$




We define $\cc$ (the union of $\{\cc(\s_\t)\}_{\t \in \ca}$ ), $\{X_\a\}_{\a \in \cc}$ and $X=\sum_{\a \in \cc} X_\a$ the same way as we did for $\ccs\ccl(n)$. Let
$$H=H_{\cc}=\sum_{\t \in \ca}\sum_{m \in Dom(\s_\t)}m(E_{\s_\t(m), \s_\t(m)}- E_{n+\s_\t(m), n+\s_\t(m)}).$$

Then $H$ is the matrix realization of $h$, hence is the  Dynkin element for the orbit $\co_{\bold d}$. And $X$ is a nilpotent element that's in $\co_{\bold d}$ (The reason is similar to the case of type $A_{n-1}$ and the reference is \cite[chap5]{CM}).

\begin{lem}\label{root}
For any roots $\a, \b \in \cc$, $\a-\b \notin \D$.
\end{lem}

Proof.  If $a$ and $\b$ are positive roots such that $\a-\b$ is a root, then we are in one of the following cases: 1) $\a=e_i\pm e_j, \b= e_i\pm e_k$ ($i \neq j, i \neq k$); 2) $\a= e_i-e_k$, $e_j-e_k$ ($j \neq l$, $i \neq k$); 3) $\a=e_i\pm e_k$ and $\b= 2e_i$ ($i \neq k$). Because of the construction of  the set $\cc$, it does not contain  two roots of the form  mentioned above.     \qed

Lemma \ref{root} is a weaker condition than Lemma \ref{chain}. Indeed, the set $\cc$ is not antichain. For example, let $\bold d=[4,2]$ and $H=diag \{ \underline 3, \overline 1, \underline 1\}$. For simplicity, we omit the negative part of the diagonal entries of $H$. Then $\cc=\{ e_1-e_3, 2e_3\} \cup \{ 2e_2\}$ and $2e_2> 2e_3$. 

Notice in the construction of the standard triple for $H$ and $X$, we only need the condition stated in Lemma \ref{root}, therefore we still can  get a standard triple $\{H, X, Y \}$ associated to $\cc$ as in Lemma \ref{trip}. Thus $\ccg_{H,i}$, $\ccq_{H, i}$,  $I_{\cc}$, $\cc^+$ and $\cc^-$ are defined accordingly as in Section 2.2.  Also we can prove that $\ccq_{H, 3}$ is contained in the ideal $I_{\cc}$ in a similar way. 

\begin{prop}
The dimension of the ideal  $I_\cc$ is equal to $m_{\co_{\bold d}}$ and $\co_{I_\cc}= \co_{\bold d}$.
\end{prop}

Proof. The main part of the  proof is basically the same as  Proposition \ref{min}. The problem is reduced to  construct a bijection between  $\cc^+$ and $\cc^-$.

Since the root space of $\ccg_{H,2}$ depends only on either  the odd parts or the even parts of $\bold d$, we  discuss  odd partitions and  even partitions separately.

Let $k_i, l_j$ be odd parts of $\bold d$ and we always have the dual indexes $k_{\tilde i_k}$ and $l_{\tilde j_l}$.  If $\a= e_{\s_{k_i}(m)}- e_{\s_{l_j}(m-2)} \notin
\cc$, where  $m >2$. Then we look at $\b=e_{k_{\tilde i_k}(m)}- e_{l_{\tilde j_l}(m-2)}$. By condition 4 of Lemma \ref{cnmap}, either $\a$ lies in $\cc^+$ and $\b$ lies in $\cc^-$ or the other way around.

By condition (5),  the two roots   $\a=e_{\s_{k_i}(2)}-e_{\s_{l_j}(0)}$ and $\b= e_{\s_{k_{\tilde i_k}}(2)}+e_{l_j(0)}$ are in bijection with each other.

Suppose $k_i, l_j$ are even parts of $\bold d$. The two roots $\a=e_{\s_{k_i}(m)}-e_{\s_{l_j}(m-2)}$ and $\b= e_{k_{\tilde i_k}(m)}-e_{l_{\tilde j_l}(m-2)}$ are in bijection with each other when $m \ge 3$ and either $i \neq \tilde i_k$ or $j \neq \tilde j_l$.

If $k_i$ and $l_j$ are even and either $i \neq \tilde i_k$ or $j \neq \tilde j_l$,then $\a=e_{\s_{k_i}(1)}+e_{\s_{l_j}(1)}$ corresponds to $\b=e_{\s_{k_{\tilde i_k}}(1)}+e_{\s_{l_{\tilde j_l}}(1)}$.

The only remaining part is $\a=e_{\s_{k_i}(m)}- e_{\s_{l_j}(m-2)}$ where $k$ and $l$ are even and $i=\tilde i_k$ and $j= \tilde j_l$. In this case, surely $k \neq l$.

If $k>l$ and $m> 5$, then  $\a$ is in bijection with $\b= e_{\s_{l_j}(m-2)}- e_{\s_{k_i}(m-4)}$
in $\D(\ccg_{H,2}) \backslash\cc$. If $m=3$, $\a$ is in bijection with $\b= e_{\s_{k_i}(1)}+ e_{\s_{l_j}(1)}$. If $k<l $ and $m \ge 3$, $\a$ corresponds to root $\b= e_{\s_{l_j}(m+2)}-e_{\s_{k_i}(m)}$.


If $i = \tilde i_k$ and $j = \tilde j_l$, then $\a=e_{\s_{k_i}(1)}+e_{\s_{l_j}(1)}$ is in bijection with $\b= e_{\s_{l_j}(3)}- e_{\s_{k_i}(1)} $, if $l>k$ or $\b= e_{\s_{k_i}(3)}- e_{\s_{l_j}(1)} $, if $l<k$.

Finally if $k$ is even,  $\a= 2e_{\s_{k_i}(1)}$ and $\b= 2e_{\s_{k_{\tilde i_k}}(1)}$ are bijective with each other.  \qed

\begin{eg}  Let $\bold d= [5^2, 3^2]$ and let $H= diag\{\underline 4, \overline 4,  \underline 2, \tilde 2, \dot 2 ,\overline 2, \underline 0, \tilde 0\}$ be the Dynkin element. Here again we omit the negative half part of $H$. $\s_{5_1}$ maps $i$ to the position of $\underline i$ and $\s_{5_2}$ maps $i$ to the position of $\overline i$ and $\s_{3_1}$ maps  $i$ to $\tilde  i$, $\s_{3_2}$ maps $i$ to $ \dot i $. \end{eg}

\section{Minimal Ideals For Type $B_n$}

Let $\ccg= \ccs \cco(2n+1)$. The set of positive roots  is $\D^+= \{e_i-e_j, e_i \mid 1 \le i<j \le n \} \cup \{e_i+e_j \mid 1 \le i,j \le n, i \neq j \} \cup \{e_i\}_{i=1}^n$. Let $\bold d$ be partition of $2n+1$ with even parts repeated with even multiplicity. As in Theorem \ref{cm1}, it corresponds to a nilpotent orbit $\co_{\bold d}$.

Suppose that $\bold d=[d^{r_d}, \dots, 2^{r_2},1^{r_1}]$, then $r_k$ is even when $k$ is an even integer. We rewrite $\bold d$ in the form  $$\bold d=[d_1, \dots, d_{r_d},\dots,  1_1, \dots, 1_{r_1}]$$ where  $\sum_{i=1}^d (r_i) i=2n+1 \text{ and } r_i \text{ is even if } 2|i.$

For each part of the partition $\bold d$, we need to attach a set of positive roots to it. Then this means again that we  should choose the appropriate index map for each $k_i$, which would determine the corresponding positive roots.
Let $\ca$ be defined as in section 2.4.

For any $k_i \in \ca$, let
\[
\s_{k_i}=\begin{cases}
   \{ k-1, k-3, \dots, 2, 0\} \rightarrow [n] \text{ if } k \text{ is odd and } 1 \le i \le \lfloor\frac{ r_{k}  }{2} \rfloor; \\
    \{ k-1, k-3, \dots, 2\} \rightarrow [n] ,\text{ if } k \text{ is odd and }\lceil \frac{ r_{k} }{2} \rceil < i \le r_{k}; \\
    \{ k-1, k-3, \dots, 1\} \rightarrow [n], \text{ if } k \text{ is even and }1 \le i \le r_{k}.
    \end{cases}
\]

In the case of $\ccs \ccp(2n)$, the odd parts have even multiplicity, therefore, we could define a dual pair of indexes $(k_i, k_{\tilde i_k})$, where $k$ is odd and $\tilde i_k:=r_k+1-i$.  The formula above defines index maps $\{\s_{k_i}\}$ for all any $k_i \in \ca$ in the last section.  However, for $\ccs \cco(2n+1)$, it's possible that $r_k$ is odd for an odd part $k$ of $\bold d$. In that case, we haven't defined the map  $\s_{k_i}$ when $i$ is equal to $\lceil \frac{ r_{k} }{2} \rceil$. Indeed, compared to the previous case, if the difficulty lies in the even parts of the partition $\bold d$ for $\ccs \ccp(2n)$, the most difficult part to construct the index map and to find the bijection for $\ccs\cco(2n+1)$ lies in its odd part.


Suppose $l^1, l^2, \dots, l^r$ are the remaining indexes  of $\ca$ without index maps associated to them.  Namely, $l^i= k_s$, where  $k$, $r_k$ are odd, and  $s=\lceil \frac{r_{k}}{2}\rceil$. Moreover, we  assume that  $l^1 < l^2 \dots < l^r$. Then $r$ must be odd since the total summation of all parts of $\bold d$ is $2n+1$. Now it's possible to define the index maps for $\{l^i \}$.


If  $i$ is odd,  set $$\s_{l^i}: \{l^i-1, l^i-3, \dots ,2 \}\rightarrow [n].$$

If $i$ is even, set $$\s_{l^i}: \{l^i-1, l^i-3, \dots , 2, 0 \}\rightarrow [n].$$

The properties for these index maps are slightly different from the previous case.

\begin{lem}\label{bnmap}
There  exits a sequence of maps $\{\s_{\t} \mid \t \in \ca\}$, satisfying the first five properties as in lemma \ref{cnmap}  with additional two:

(6) Let $k,l$ be odd integers, $1 \le i \le \lfloor r_k/2 \rfloor$ and $j= \tilde j_l$, then $\s_{k_i}(m)<\s_{l_j}(m)<\s_{k_{\tilde i_k}(m)}$. If $m=0$, then $\s_{k_i}(0)< \s_{l_j}(0)$.

(7) If $k,l $ are odd integers, $k<l$ and $i=\tilde i_k$, $j=\tilde j_l$, then $\s_{k_i}(m)<\s_{l_j}(m)$.

\end{lem}

The last two properties give additional restrictions for the placements of the indexes that come from the odd parts of $\bold d$ with odd multiplicities.

For the even part $k_i$ of $\bold d$,   if $1 \le i \le \frac {r_k}{2}$,  we can attach a set of  positive roots to the map $\s_{k_i}$:
\begin{align*} \cc(\s_{k_i})= \{e_{\s_{k_i}(k-1)}- e_{\s_{k_i}(k-3)} , \dots, e_{\s_{k_i}(3)}- e_{\s_{k_{i}}(1)}, e_{\s_{k_i}(1)}+ e_{\s_{k_{\tilde i_k}}(1)}. \} \end{align*}

If $k$ is even and $ i > \frac{r_k}{2}$,  set:

\begin{align*}\cc(\s_{k_i})= \{ e_{\s_{k_i}(k-1)}- e_{\s_{k_i}(k-3)} , \dots, e_{\s_{k_i}(3)}- e_{\s_{k_i}(1)} \}.
\end{align*}


For any odd part $k_i$ of $\bold d$, if $i < \tilde i_k$, set
\begin{align*}
\cc(\s_{k_i})= \{ e_{\s_{k_i}(k-1)}-e_{\s_{k_i}(k-3) }, \dots, e_{\s_{k_i}(2)}- e_{\s_{k_i}(0)}\}.
\end{align*}

If $i> \tilde i_k$, the set of positive roots we attached to the map $\s_{k_i}$ is
\begin{align*}
\cc(\s_{k_i})= \{ e_{\s_{k_i}(k-1)}-e_{\s_{k_i}(k-3) }, \dots, e_{\s_{k_i}(2)}+ e_{\s_{k_{\tilde i_k}}(0)}\}.
\end{align*}

The remaining case is that there exists some $i$ such that $i= \tilde i_k$. In this case, if $0 \in Dom(\s_{k_i})$, we can attach
\begin{align*}
\cc(\s_{k_i})= \{ e_{\s_{k_i}(k-1)}-e_{\s_{k_i}(k-3) }, \dots,  e_{\s_{k_i}(2)}+ e_{\s_{k_i}(0)}, e_{\s_{k_i}(2)}- e_{\s_{k_i}(0)}\}
\end{align*}
to the map $\s_{k_i}$. Otherwise, $0 \notin Dom(\s_{k_i}) $,  we can attach almost the same chunk of roots to $\s_{k_i}$ except replacing the last two roots in $\cc(\s_{k_i})$ with $e_{\s_{k_i}(2)}$. Namely \begin{align*}
\cc(\s_{k_i})= \{ e_{\s_{k_i}(k-1)}-e_{\s_{k_i}(k-3)}, \dots, e_{\s_{k_i}(4)}-e_{\s_{k_i}(2)}, e_{\s_{k_i}(2)}\}
\end{align*}

For each map $\s_{k_i}$, once we get the set of positive roots associated to $\s_{k_i}$, we can define $H, X, \cc, I_{\cc}, \ccg_{H, i}, \ccq_{H, i}, \cc^+ $ and $\cc^-$ the same way as in previous section.  The fact that $\co_X=\co_{\bold d}$ comes from  \cite[Chap 6]{CM}. In the case of type $B_n$, $\cc$ is not an antichain and does not satisfy lemma \ref{root}.  It's because that  $e_i$, $e_j$ can be both in $\cc$ for some $i, j$. That means  we cannot get a standard triple from $H, X$.  However,  we will prove

\begin{prop}\label{dim}The ideal $I_{\cc}$ has associated orbit $\co_{\bold d}$ and $\dim I_{\cc}= m_{\co_{\bold d}}$.
\end{prop}

This follows from Proposition \ref{253} below.

\begin{prop}\label{253}
There exists a bijection between $\cc^+$ and $\cc^-$.
\end{prop}

Proof. We discuss the even parts and odd parts of the partition $\bold d $ separately. For simplicity, we will use $k_i$ to denote the map $\s_{k_i}$.

Suppose $k_i, l_j$ are even parts of the partition $\bold d$ and $k_i \neq l_j$ as index in $\ca$.  Then  $\a=e_{{k_i}(m)}-e_{{l_j}(m-2)}$ is in bijection with $\b=e_{{k_{\tilde i_k}}(m)}-e_{{l_{\tilde j_l}}(m-2)}$ when $m>2$. Namely either $\a \in \cc^+ $ and $\b \in \cc^-$ or $\a \in \cc^-$ and $\b \in \cc^-$ (by condition 4 of lemma \ref{bnmap}). If $\a= e_{{k_i}(1)}+e_{{l_j}(1)}$, then it is in bijection with $\b= e_{{k_{\tilde i_k}}(1)}+e_{{l_{\tilde j_l}}(1)}$ (also by lemma \ref{bnmap}).

Suppose that $k_i, l_j$ are odd parts of $\bold d$. If $m >2 $ and either $i \neq \tilde i_k$ or $j \neq \tilde j_l$, then the roots $\a=e_{{k_i}(m)}-e_{{l_j}(m-2)}$  and $\b=e_{{k_{\tilde i_k}}(m)}-e_{{l_{\tilde j_l}}(m-2)} $ are bijective with each other. If $\a= e_{{k_i}(2)}$, where $i \neq \tilde i_k$, the root  corresponding to $\a$ is $\b= e_{k_{\tilde i_k}(2)}$. If $\a= e_{{k_i}(2)}-e_{\s_{l_j}(0)}$, then $\b = e_{\s_{k_{\tilde i_k}}(2)}+e_{\s_{l_j}(0)}$ is the corresponding root.

The final step is to find the bijection when $k \neq l$, $i=\tilde i_k$ and $j= \tilde j_l$. Suppose that  $m > 4$ and $k>l$, the root $\a=e_{k_i(m)}-e_{l_j(m-2)}$ lies in $\cc^-$ while its corresponding root $\b= e_{l_j(m-2)}-e_{k_i(m-4)}$ lies in $\cc^+$. For any $\a= e_{k_i(2)}-e_{l_j(0)}$ that lies in $\cc^-$, i.e. $k>l$, it corresponds to a positive root $\b= e_{k_i(2)}+e_{l_j(0)}$ that lies in $\cc^+$.

The remaining positive root that lies in the set $\cc^-$ has the form  $\a= e_{k_i(4)}-e_{l_j(2)}$, with $i=\tilde i_k$ and $j= \tilde j_l$. Then we need to switch our notation to $\a=e_{l^s(4)}-e_{l^t(2)}$, where $l^s=k_i$ and $l^t=l_j$. Since $\a \in \cc^-$, by condition 7 of lemma \ref{bnmap}, we have $l^s>l^t$ and $s>t$.

In this case, if $0 \notin Dom(l^t)$, then $e_{l^t(2)} \in \cc$, and $\a$ corresponds to $\b= e_{l^t(2)}- e_{l^s(0)} \in \cc^+$ if $0 \in Dom(l^s)$ or $\a$ corresponds to $\b=e_{l^t(2)}- e_{l^{s-1}(0)} \in \cc^+$ if $0 \notin Dom(l^s)$.

If $0 \in Dom(l^t)$, then $e_{l^t(2)}\pm e_{l^t(0)} \in \cc$, and $\a$ corresponds to $\b= e_{l^t(2)}- e_{l^s(0)}$ if $0 \in Dom(l^s)$ or $\a$ corresponds to $\b=e_{l^t(2)}- e_{l^{s-1}(0)}$ if $0 \notin Dom(l^s)$ and $s> t+1$. If $s=t+1$, then $0 \notin Dom(l^s)$ and $a$ corresponds to $\b= e_{l^t(2)} \in \cc^+$. \qed



\section{Minimal Ideals For Type $D_n$}

Let $\ccg= \ccs \cco(2n)$. The matrix realization of $\ccs\cco(2n)$ is
\[\{ \left( \begin{array}{cc}
 Z_1 & Z_2 \\
  Z_3 & -Z_1^t
 \end{array}\right) \mid Z_i \in M_{n}(\mathbb C), Z_2, Z_3 \text{ skew-symmetric } \}.
\]

The set of positive roots  is $\D^+= \{e_i-e_j\mid 1 \le i<j \le n \} \cup \{e_i+e_j \mid 1 \le i,j \le n, i \neq j \}$. Let $\bold d$ be partition of $2n$ with even parts repeated with even multiplicity. If $\bold d$ is not a very even partition,  as in Theorem \ref{cm1}, it corresponds to a nilpotent orbit $\co_{\bold d}$. If $\bold d $ is a very even partition, then it corresponds to two nilpotent orbits $\co^I_{\bold d}$ and $\co^{II}_{\bold d}$.

Suppose that $\bold d=[d^{r_d}, \dots, 2^{r_2},1^{r_1}]$, then $r_k$ is even when $k \in 2\mathbb N$. We rewrite $\bold d$ in the form  $$\bold d=[d_1, \dots, d_{r_d},\dots, 1_1, \dots, 1_{r_1}], \text{ where } \sum_{i=1}^d (r_i)i=2n \text{ and } r_i \text{ is even if } 2|i.$$

Recall the procedure to get the weighted Dynkin diagram from the  partition $\bold d$. We can obtain a sequence of integers from $\bold d$ the same way as we did for $\ccg=\ccs\ccp(2n)$.  The sequence takes the form $$(h_1, \dots, h_n, -h_1, \dots, -h_n), \text{ where } h_1 \ge h_2 \ge \dots \ge h_n. $$ If $\bold d$ is not a very even partition, following previous matrix realization of $\ccg$, the Dynkin element for the orbit $\co_{\bold d }$ is $$H= diag\{h_1, h_2, \dots, h_n, -h_1, \dots, -h_n \}.$$  However, if $\bold d$ is a very even partition, the Dynkin elements for the orbits $\co^{I}_{\bold d}$ and $\co^{II}_{\bold d}$ are $$H_1=diag\{h_1,  \dots, h_n, -h_1, \dots, -h_n \}$$ and $$H_2 =diag\{h_1,  \dots, h_{n-1}, -h_n, -h_1, \dots, -h_{n-1},  h_n \},$$  both of which are dominant.

Let $\ca$ be the same as in last section and we define the index maps $\{\s_{k_i}\}_{k_i \in \ca}$ the same way as we did for $\ccg=\ccs\cco(2n+1)$. The notations $l^1, l^2, \dots, l^r$ have the same meaning as in last section. The only difference is that  $r$ is even. 

\begin{lem}
There  exits a sequence of maps $\{\s_{\t} \mid \t \in \ca\}$, satisfying the first five properties as in lemma \ref{cnmap}  with additional two:

(6) Let $k,l$ be odd integers, $1 \le i \le \lfloor r_k/2 \rfloor$ and $j= \tilde j_l$, then $\s_{k_i}(m)<\s_{l_j}(m)<\s_{k_{\tilde i_k}(m)}$. If $m=0$, then $\s_{k_i}(0)< \s_{l_j}(0)$.

(7) If $k,l $ are odd integers, $k<l$ and $i=\tilde i_k$, $j=\tilde j_l$, then $\s_{k_i}(m)<\s_{l_j}(m)$.

\end{lem}

For every even part $k_i$ of $\bold d$ and for  odd part $k_i$ such that $i \neq \tilde i_k:=r_k-i+1$, we attach the same set of positive roots $\cc(\s_{k_i})$ to $\s_{k_i}$ as in section 2.5.

For $l^1 \le \dots \le l^r$, there is always a pair $(l^i, l^{i+1})$, where $1 \le i \le r$ and $i$ is odd. 

If $i$ is odd, we attach the following set of roots to $\s_{l^i}$:
\begin{align*}
\cc(\s_{l^i})= \{ e_{\s_{l^i}(l^i-1)}-e_{\s_{l^i}(l^i-3) }, \dots, e_{\s_{l^i}(2)} \pm  e_{\s_{l^{i+1}}(0)}\}.
\end{align*}

If $i$ is even, we attach the following set of roots to $\s_{l^i}$:
\begin{align*}
\cc(\s_{l_i})= \{e_{\s_{l^i}(l^i-1)}-e_{\s_{l^i}(l^i-3) }, \dots, e_{\s_{l^i}(2)} - e_{\s_{l^i}(0)}\}.
\end{align*}

Similarly, we have $\cc, X, H, I_\cc, \cc^+, \cc^-$ and $I_\cc$. If $\bold d$ is not an even partition, $\co_X$ is the unique orbit $\co_{\bold d}$.

If $\bold d$ is a very even partition, then $X$ corresponds to  the first orbit $\co_{\bold d}^I$ and the Dynkin element $H$ has the form $$H= diag\{h_1,  \dots, h_n, -h_1, \dots, -h_n \}.$$ To get a representative for another orbit, first notice that  $n=\s_{k_i}(1)$ for some $k_i \in \ca$ and $k$ is even. Let $j=\tilde i_k=r_k-i+1$. We let  \begin{align*}
 \tilde \cc(\s_{k_i})\cup \tilde \cc( \s_{k_j}) =&
\{ e_{\s_{k_i}(k-1)}- e_{\s_{k_i}(k-3)} , \dots, e_{\s_{k_i}(3)}+ e_{\s_{k_i}(1)}, e_{\s_{k_i}(1)}- e_{\s_{k_j}(1)} \} \\& \cup \{  e_{\s_{k_j}(k-1)}- e_{\s_{k_j}(k-3)} , \dots, e_{\s_{k_j}(3)}- e_{\s_{k_j}(1)}\}.
\end{align*}

For any  $\t \in \ca$, if $\t\neq k_i, k_j$, let $\tilde \cc(\s_\t)=\cc(\s_\t)$ and $\tilde \cc= \cup_{\t \in \ca}\tilde \cc(
s_\t)$.  Let $\tilde X= \sum_{\t \in \tilde \cc}X_\t$, it corresponds to the orbit $\co_{\bold d}^{II}$. The Dynkin element $\tilde H$ has the form $\tilde H=diag\{h_1,  \dots, h_{n-1}, -h_n, -h_1, \dots, -h_{n-1},  h_n \}$.

In both cases, it's easy to see that $\cc$ and $\tilde \cc$ are antichains and there exist triples $\{H, X ,Y\}$ and $\{\tilde H, \tilde X, \tilde Y\}$.

\begin{eg} Let $\bold d=[4^2]$. Then $H=\diag(\overline 3, \underline 3, \overline 1, \underline 1, -3,-3, -1,-1)$, where $\{\overline 3, \overline 1\}= Im(\s_{4_1})$ and $\{\underline 3, \underline 1 \}=Im(\s_{4_2})$. Then $\cc= \{e_1-e_3, e_3+e_4, e_2+e_4\}$ and $\tilde \cc=\{e_1-e_3, e_3-e_4, e_2+e_4\}$.
\end{eg}
We can obtain Proposition \ref{dim} in a similar way for $\ccs\cco(2n)$. If $\bold d$ is a very even partition, then $I_\cc$ is the minimal ideal for $\co_{\bold d}^I$ and $I_{\tilde \cc}$ is the minimal ideal for $\co_{\bold d}^{II}$. Here we omit the proof.

\bibliographystyle{amsalpha}

\end{document}